\documentclass{psapm-l}
\usepackage{amsmath, amsfonts, amsthm, amssymb, verbatim} 
\newcommand \lambdat {\widetilde \lambda}
\newcommand \lambdab {\overline \lambda}
\newcommand \at {\widetilde a}
\newcommand \lamt {\widetilde \lambda}
\newcommand \lt {\widetilde l}
\newcommand \rt {\widetilde r}
\newcommand \sgn {\text{sgn}}
      
\newcommand \del   {{\partial}} 
\newcommand \eps   {\varepsilon}     
\newcommand \RR   	{{\mathbb R}}   
\newcommand \R   	{{\mathbb R}}  
\newcommand \ab {\overline a}
\newcommand \be {\begin{equation}}
\newcommand \ee {\end{equation}}   
\newcommand \dist {\text{dist}}

\newcommand \At {\widetilde A}

\newcommand \RN    {\mathbb{R}^N}

\newcommand \Mcal {\mathcal M}
\newcommand \Lcal {\mathcal L}
\newcommand \NN {\mathcal N}  
   
\newcommand \Dbold {\mathbf D}   
\newcommand \Abar  {{\overline A}}
\newcommand \abar  {{\overline a}}
\newcommand \Ccal {\mathcal C}
\newcommand \Jcal {\mathcal J}
\newcommand \Ical {\mathcal I} 
\newcommand \Scal {\mathcal S}
\newcommand \Fcal {\mathcal F}
\newcommand \Rcal {\mathcal R}
\newcommand \lam \lambda 
\newcommand \lamb {\overline \lambda}
\newtheorem{theorem}{Theorem}[section]
\newtheorem{lemma}[theorem]{Lemma}
\theoremstyle{definition}
\newtheorem{definition}[theorem]{Definition}
\newtheorem{example}[theorem]{Example}
\newtheorem{proposition}[theorem]{Proposition}

\theoremstyle{remark}

\numberwithin{equation}{section}

\begin{document}

\title[Stability in the $L^1$ norm via a linearization method]
{Stability in the $L^1$ norm via a linearization method for nonlinear hyperbolic systems}
 
\author[P.G. L{\tiny e}Floch]{Philippe G. LeFloch}

\address{Laboratoire Jacques-Louis Lions \& Centre National de la Recherche Scientifique (CNRS),
Universit\'e Pierre et Marie Curie (Paris 6), 4 Place Jussieu,  75252 Paris, France.} 
\email{pgLeFloch@gmail.com} 

\date{November 2008.}
 

\begin{abstract}
We discuss the existence and uniqueness of discontinuous solutions to adjoint problems 
associated with nonlinear hyperbolic systems of conservation laws. By generalizing the Haar method 
for Glimm-type approximations to hyperbolic systems, we establish that entropy solutions 
depend continuously upon their initial data in the natural $L^1$ norm. 
\end{abstract}

\maketitle


\section{Introduction}
  
We discuss the existence and uniqueness of discontinuous solutions to adjoint problems associated 
with nonlinear hyperbolic systems of conservation laws. By generalizing the Haar method \cite{Haar}
to encompass discontinuous solutions generated by Glimm \cite{Glimm} and
front-tracking \cite{Dafermos1, DiPerna, Risebro, Bressan1} approximation methods, 
we establish that entropy solutions depend continuously upon their initial data in the natural $L^1$ norm. 
In order to implement the proposed {\sl linearization method,}  
we analyze linear hyperbolic systems with discontinuous coefficients, possibly in a 
non-conservative form. Our analysis begins with the key observation 
that while entropy solutions, by definition, contain compressive shocks only, the 
averaged matrix 
$$
\Abar(u, v) = 
\int_0^1 Df (u + \theta(v - u)) \, d\theta
$$
(associated with the flux of the system $f$ and two entropy solutions $u, v$) may contain
compressive or undercompressive shocks, but no rarefaction-shocks. 
This is an essential observation since, otherwise, 
rarefaction shocks would be a source of non-uniqueness and instability. 
The proposed method rests on geometric properties of $\Abar(u, v)$ and on the 
construction of a weighted norm which is determined during the evolution and 
takes into account wave cancellation effects along generalized characteristics.
This strategy was presented by the author in 1998 (in a lecture  
at the University of Wisconsin-Madison) and appeared in print 
in \cite{HuLeFloch,LeFloch-book}. 
In independent work, 
another proof of this continuous dependence property for genuinely nonlinear systems, 
based on earlier work on scalar conservation laws by Liu and Yang \cite{LiuYang1},  
was obtained simultaneously to \cite{HuLeFloch,LeFloch-book} 
by Bressan, Liu, and Yang \cite{BressanLiuYang,LiuYang2}. 
In the present review, we cover the recent results on non-genuinely nonlinear systems established 
in \cite{LeFloch-cancel,LeFloch-prepare}.


\section{Continuous dependence property via a linearization method}
\label{section22}

We consider solutions with small total variation, constructed by Glimm or front tracking schemes, 
satisfying the strictly hyperbolic system 
\be
\label{11}
\del_t u + \del_x f(u) = 0, \quad u=u(t,x) \in \RN.  
\ee
 The flux $f(u) \in \RN$ need not be genuinely nonlinear and we solely assume that, for all relevant values
 of $u$, the matrix  
 $Df(u)$ 
has distinct eigenvalues $\lambda_j(u)$ and basis of left- and right-eigenvectors $l_j(u)$ and $r_j(u)$, respectively. 
Given two entropy solutions $ u,v$ we define 
$$
\psi := v -u
$$
and introduce the averaged matrix
$$
\At := \Abar(u,v) := \int_0^1 Df(\theta u + (1-\theta) v) \, d\theta.
$$
Clearly, if $u,v$ satisfy \eqref{11}, then $\psi$ is a solution to the following linear hyperbolic system
\be
\label{12}
\del_t \psi + \del_x ( \At \, \psi) = 0. 
\ee

We observe that the $L^1$ stability property for the linear system, i.e. 
$$ 
\| \psi(t) \|_{L^1(\RR)}  \lesssim \|\psi(0)\|_{L^1(\RR)}
$$
for a sufficiently large class of matrices $ \At$ and solutions $ \psi$, 
implies the $L^1$ continuous dependence property for the nonlinear system 
\be
\label{Lone}
\| (u - v)(t) \|_{L^1(\RR)}  \lesssim \|(u - v)(0)\|_{L^1(\RR)}. 
\ee 
The main difficulties are, first, to identify a suitable class of linear systems and solutions
and, second, to construct a weighted $L^1$ norm that decreases in time.  
Our objective in this short presentation is to briefly review our linearization method 
and present several new results; we refer to \cite{LeFloch-cancel,LeFloch-prepare} for further details.


Throughout, it is convenient to work with piecewise constant data and solutions. We suppose that  
the given matrix-valued field 
$\At$ is strictly hyperbolic with distinct eigenvalues $\lambdat_j$ and eigenvectors $\lt_j, \rt_j$. 
The points of continuity, jump, and interaction in $\At$ are denoted by 
$$
\RR_+ \times \RR = \Ccal(\At) \cup \Jcal(\At) \cup \Ical(\At),
$$
respectively.  

\begin{definition} \it 
\label{sys}
A $j$-discontinuity $(t,x) \in \Jcal(\At)$ propagating at the speed $\lamb$ 
is said to be 
\begin{enumerate}

\item[--] {\rm compressive} if $ \lambdat_j^- \geq \lamb \geq \lambdat_j^+$,

\item[--] {\rm slow undercompressive} if $\lamb < \min\bigl(\lambdat_j^-, \lambdat_j^+\bigr)$, 
 
\item[--] {\rm fast undercompressive} if $ \lamb > \max\bigl(\lambdat_j^-, \lambdat_j^+\bigr)$, or 
 
\item[--] {\rm a rarefaction-shock}  $ \lambdat_j^- \leq \lamb \leq \lambdat_j^+$. 
\end{enumerate}

\end{definition}

We use the following decomposition 
$$
 \Jcal(\At) =: \Lcal(\At) \cup \Scal(\At) \cup \Fcal(\At) \cup \Rcal(\At).
$$
We first point out that if $\At$ contains rarefaction-shocks then
the uniqueness property and, consequently,
the continuous dependence property fail, as illustrated by the following example.

\begin{example} \it 
Consider the case that the speed coefficient is a rarefaction shock in a scalar equation, i.e. 
$$ 
\del_t \psi + \del_x ( \at \, \psi) = 0,  \quad \psi=\psi(t,x) \in \RR, 
$$
with   
$$ 
\at(t,x) = \begin{cases}
 -1,  &   x<0, 
\\
\hskip.23cm 1,  &   x>0. 
\end{cases}
$$
Then, the corresponding Cauchy problem admits infinitely many weak solutions, since 
within the wedge $|x/t| < 1$ the solution
can not uniquely determined from the given initial data by the method of characteristics.  
\end{example}

\

Let us summarize key facts of the proposed linearization method:  

\begin{enumerate}

\item Rarefaction shocks may lead to instability and non-uniqueness. 

\item However, the entropy condition assumed by the weak solutions $u,v$ to nonlinear hyperbolic 
systems implies that all shock waves in $u,v$ are {\sl compressive}
and, more importantly, that the averaged matrix $\Abar(u,v)$ can not contain rarefaction shocks.   

\item $L^1$ stability estimates for rarefaction-free systems  
are established by defining functionals that are equivalent to the $L^1$ distance 
and are generated along the time-evolution by a constructive approach.  

\item Furthermore, in the course of this analysis, pointwise convergence properties of Glimm-type schemes
are required to establish the stability of certain nonconservative products.  

\end{enumerate}

Details can be found in \cite{HuLeFloch,LeFloch-book,GoatinLeFloch1,GoatinLeFloch2,LeFloch-cancel,LeFloch-prepare}. 
This method was first investigated, in the $L^2$ norm rather than the $L^1$ one, 
in the earlier papers \cite{Oleinik,LeFlochXin}. 

Let us conclude this section by recalling the pointwise convergence property 
of Glimm-type methods, as established by Glimm and Lax \cite{GlimmLax} and DiPerna \cite{DiPerna,DiPerna1}.  
See also Dafermos \cite{Dafermos2,Dafermos-book}, Liu \cite{Liu2}, 
and  Bressan and LeFloch \cite{BressanLeFloch2}. 

When applied to the nonlinear hyperbolic system \eqref{11} with initial data of small total variation,
Glimm solutions, say $u^h=u^h(t,x)$, {\sl converge locally uniformly} to the entropy solution $u$
of the corresponding Cauchy problem, 
in the 
sense that the following two properties hold for all but countably many $t_0$ and all $ x_0$:  
\begin{enumerate}
\item
         If $(t_0, x_0)$ continuity point for $u$, then for all $ \eps>0$ 
         there exists a neighborhood $\NN(t_0, x_0)$ of that point and a real $h_0>0$ such that 
for all $h < h_0^\eps$ and $(t,x) \in \NN^\eps(t_0, x_0)$ 
$$ 
|u^h(t,x) - u(t_0,x_0)| + |u(t,x) - u(t_0,x_0)| < \eps.
$$

\item  
        If $ (t_0, x_0)$ a discontinuity point for $u$, then  there exists
         a shock curve $ t \mapsto y(t)$ defined near $t_0$ with $y(t_0) = x_0$ 
and, for $ \eps>0$, there exist a neighborhood $ \NN^\eps(x_0, t_0)$,
 a real $h_0^\eps>0$, and an approximate shock curve $ t \mapsto y^h(t)$ 
such that for $h<h_0^\eps$ and $(t,x) \in \NN^\eps(t_0, x_0)$ 
$$
|y^h(t) - y(t)| < \eps,  
$$ 
provided  $(t,y^h(t)), \, (t,y(t)) \in \NN^\eps(t_0, x_0)$, 
and for $(t,x) \in \NN^\eps(t_0, x_0)$  
$$
\aligned
|u^h(t,x) - u(t_0,x_0\pm)|  &< \eps,    \quad  && x \gtrless y^h(t),  
\\
  |u(t,x) - u(t_0, x_0\pm)|  &< \eps,    \quad  && x \gtrless y(t).
\endaligned   
$$ 
\end{enumerate}

\section{Linear and nonlinear scalar equations} 

\subsection*{Linearization involving one solution}
 
Consider the two scalar equations
\be
\label{eqq}
\del_t \psi + \del_x ( \widehat{f_u}(u) \, \psi) = 0, \qquad \del_t u + \del_x f(u)= 0, 
\ee
where $f:\RR \to \RR$ is strictly convex.
We emphasize that, when $ u=v$, all shocks in $\at=\abar(u,u)= f_u(u)$ are {\sl compressive.}

Following \cite{LeFloch-IMA1}, weak solutions $\psi$ are sought as measures in the spatial variables 
and \eqref{eqq} is defined in the sense of Volpert \cite{Volpert}, which is equivalent to choosing the 
family of straightlines in Dal~Maso-LeFloch-Murat's theory \cite{DLM}.  
We impose initial data $\psi_0, u_0$ such that $u_0 \in BV(\RR)$, the space of functions with 
bounded variation, and $\psi_0 \in \Mcal_b(\RR)$, the space of bounded measures.   

First, we consider the Riemann problem associated with data $u_l, u_r$ and $ \psi_l, \psi_r$, respectively.  
We distinguish between two cases: either $u$ is a shock (with speed denoted by $\lambdab$)
or $u$ is a rarefaction.

\begin{theorem}[Riemann problem \cite{LeFloch-IMA1}]
\label{bof} 
\begin{enumerate}
\item If the solution $u$ is a shock, then the solution $\psi$ to the corresponding Riemann problem
is given by 
$$
\psi(t,x) = \psi_l + (\psi_r - \psi_l) H(x-\lamb t) + t \, (C_r - C_l) \delta_{x-\lamb t},  
$$
with 
$$ 
C_r := (\lamb - f_u(u_r)) \psi_r, \qquad C_l:= (\lamb - f_u(u_l)) \psi_l.
$$ 

\item If the solution $u$ is a rarefaction, the Riemann problem admits infinitely many solutions, 
in particular ($ \varphi_* \in \RR)$  
$$
\aligned
 \psi(t,x) = 
\, &  \psi_l \, (1- H(x-t f_u(u_l)) + \psi_r H(x-t f_u(u_r))  
\\
&   + \varphi_* \, \delta({x-t f_u(u_r)}) - \varphi_*  \delta({x-t f_u(u_l)}). 
\endaligned
$$
\end{enumerate}

\end{theorem}

The following remarks are in order:  
 
\begin{enumerate} 

\item Case of a shock wave: 
\begin{enumerate} 
\item[--] The Riemann problem admits a unique solution. 

\item[--] A Dirac mass in $\psi$ propagates along the shock trajectory and its
 strength grows linearly in time, and the solution is solely a bounded measure. 

\item[--] However, under the compatibility condition $ C_l = C_r$, the solution $\psi$ 
is more regular and is a function of bounded variation. 

\end{enumerate}

\item Case of a rarefaction wave: 
\begin{enumerate} 

\item[--] The Riemann problem admits (at least) a one-parameter family of solutions. 

\item[--] The solutions $\psi$ contain two propagating Dirac masses and are bounded
measures, only. 

\item[--] Still, the $\psi$-equation does admit bounded variation solutions, 
obtained for instance by taking $\varphi_*=0$ in Theorem~\ref{bof}.  

\end{enumerate}
\end{enumerate}


 We now turn our attention to the Cauchy problem and we search 
 for solution that are bounded measures in the space variable. 

\begin{theorem}[Existence and uniqueness theory \cite{LeFloch-IMA1}]
\label{56}
\begin{enumerate}
Let $u$ be an entropy solution with bounded variation to a conservation law with convex flux $f$. 

\item For every initial data $\psi_0 \in L^1(\RR) \cap L^\infty(\RR)$, the initial value problem for  
the equation 
$$
\del_t \psi + \del_x ( \widehat{f_u}(u) \, \psi) = 0
$$
admits 
at least one solution $\psi \in L^\infty(\RR_+, \Mcal_b(\RR))$, 
provided the product $\widehat{f_u}(u) \, \psi$ is understood as a Volpert product. 

\item When $\psi_0 \in \Mcal_b(\RR)$  and $\del_x u_0 \leq C_0$, the problem admits at most one such solution. 

\end{enumerate}

\end{theorem}

The proof of the existence part uses Lax formula for the entropy solution $u$
and allows us to write an explicit formula for $\psi$ by tracking forward the initial data $\psi_0$ along 
the generalized characteristics associated with $u$.  
The uniqueness part is based on an $L^1$-type contraction argument, which strongly uses
the entropy condition satisfied at jump discontinuities.

A large literature is available on linear hyperbolic equations. 
More general speed coefficients are now covered by Bouchut and James \cite{BouchutJames}
and Popov \cite{Popov}. The connection with the study of finite difference schemes 
was extensively investigated by Tadmor \cite{Tadmor}.  
Important progress was also made on multidimensional equations 
by Ambrosio \cite{Ambrosio}, Colombini, Crippa, and Rauch \cite{CCR}; 
these latter papers impose a restriction on the divergence of the velocity field 
and do not cover the equations and solutions in Theorem~\ref{56}.

  
\subsection*{Linearization involving two solutions}

Dealing with two solutions is more delicate, since now the average speed does not contain only  
compressive shocks. 
 
\begin{theorem}[Sharp $L^1$ stability property \cite{LeFloch-prepare}]
Consider the linear hyperbolic equation 
$$
\del_t \psi + \del_x (\at \, \psi) = 0,
$$ 
where $ \at :=\abar(u,v)$ and $ u,v$ are bounded variation entropy solutions 
of a scalar conservation law with general flux $f:\RR \to \RR$. 
Then, for all BV solutions $\psi$ generated by front tracking one has 
$$  
\aligned
& \|\psi(t)\|_{L^1(\RR)} +   \Dbold_2(t) + \Dbold_3(t)
\lesssim \|\psi(0)\|_{L^1(\RR)}
\\
& \Dbold_2(t) :=  \int_0^t \sum_{\Lcal(\at)}  | \lam^{\at} - \at_-| \, |\psi_-| \, d\tau. 
\\
& \Dbold_3(t) 
 := \sum_{\Scal(\at) \cup \Fcal(\at)} 
 | \at_- - \lamb |\, |\at_+ - \at_-| \, |\psi_-| + \int_\RR |\at - f'(u) | \, |\psi| \, dV_{\at}^c. 
\endaligned 
$$

\end{theorem} 

Here, $V_a^c$ denotes the absolutely continuous part of the total variation measure of $a$. 
It should be observed that the term $\Dbold_2(t)$ provides the strongest, quadratic decay and is
associated with Lax shocks. The term $\Dbold_3(t)$ is cubic in nature and
involves undercompressive waves $\Scal(\at) \cup \Fcal(\at)$, only.   
 
The above result immediately applies to $\psi = v -u$ and yields a {\sl sharp continuous dependence property}
which generalizes the standard Kruzkov's contraction property.

Earlier results about the continuous dependence property 
for scalar conservation laws concerned {\sl convex flux,} only. 
\begin{enumerate}

\item For approximate piecewise smooth solutions constructed by the Glimm scheme, 
Liu and Yang \cite{LiuYang1} introduced an explicit functional in order to control the $L^1$ norm between two solutions.  

\item For exact solutions with bounded variation, Dafermos \cite{Dafermos-book}
first derived a sharp version using the method of generalized characteristics. 
Later, Goatin and LeFloch \cite{GoatinLeFloch1} covered a class of hyperbolic equations and developed the technique
based on compressive and undercompressive discontinuities.  

\end{enumerate}

In \cite{LeFloch-prepare}, we show that, 
for scalar conservation laws with arbitrary flux-function, the averaged coefficient $\at$ contains no rarefaction shocks. 
We construct a decreasing weighted $L^1$ norm for the solution $\psi$.  We establish that 
the nonconservative products 
(especially the term in $ \Dbold_3$)
arising in the sharp stability estimate
are stable,  
thanks to the pointwise convergence property of the Glimm-type schemes mentioned earlier.  
Such argument were first developed by LeFloch and Liu \cite{LeFlochLiu} in their
 version of the Glimm scheme for nonconservative systems. (See Section~\ref{DLMsection}.)

For scalar conservation laws, Definition~\ref{DLMsection} simplifies into the following form.

\begin{definition} \it
A shock wave at a point $ (t,x) \in \Jcal(\at)$ associated with the speed $ \ab$
is said to be  
\begin{enumerate}
\item[--] {\rm compressive} if $ \at_- \geq \ab \geq \at^+$, 

\item[--] {\rm slow undercompressive} if $ \ab < \min\bigl(\at^-, \at^+\bigr)$,  
 
\item[--] {\rm fast undercompressive}  if  $ \ab > \max\bigl(\at^-, \at^+\bigr)$, or 

\item[--] {\rm a rarefaction-shock}  if $ \at^- \leq \ab \leq \at^+$. 
\end{enumerate} 
\end{definition}

The following classification provides us with a characterization of the nature of shocks
in the averaged speed coefficient associated with two entropy solutions.
We consider a scalar conservation law with flux $ f:\RR \to \RR$ and a propagating
discontinuity $(u_-, u_+)$ whose propagation speed is 
$$ 
\abar(u_-, u_+) := \int_0^1 \del_u f(s \, u_- + (1- s) \, u_+) \, ds. 
$$
Fix a constant $v \in \RR$, representing the (constant) local value of another solution. 
Then, whether the discontinuity   
$$
 \at(t,x) = \begin{cases} 
 \abar_- := \abar(u_-,v),  &   x <  \abar(u_-, u_+) \, t, 
\\
 \abar_+ := \abar(u_+,v),  &   x >  \abar(u_-, u_+) \, t, 
\end{cases}
$$
is a compressive, undercompressive, or rarefaction shock 
is uniquely determined by the sign of 
$$
 \Omega := (u_\pm - v) \, (u_+ - u_-) \, (\abar_+ - \abar_-).
$$ 
Specifically, we have: 
\begin{enumerate} 
\item[$\bullet$] If $ u_\pm - v$ have the same sign, then 
$$ 
\Omega = (u_- - v) \, (u_+ - u_-) \, (\ab_+ - \ab_-) 
\begin{cases} 
 \leq 0, &  \Scal(\at),     
\\ 
 \geq 0, &  \Fcal(\at). 
\end{cases}  
$$
\item[$\bullet$] If $ u_\pm - v$ have opposite signs, then 
$$ 
\Omega = (u_- - v) \, (u_+ - u_-) \, (\ab_+ - \ab_-) 
\begin{cases} 
 \geq 0, &  \Lcal(\at),  
\\
 \leq 0, &  \Rcal(\at).  
\end{cases}  
$$
\end{enumerate}

Next, by recalling that solutions under consideration also satisfy the entropy condition, 
we arrive at: 

\begin{proposition}[Fundamental property for scalar equations]
Let~$ u,v$ be two entropy solutions (with arbitrarily large total variation) to 
a scalar conservation law with general flux.   
Then, an (entropy admissible) shock wave of $u$ or $v$ {\rm can not be a rarefaction-shock} of 
the averaged speed $ \at := \abar(u,v)$. 
\end{proposition}

\section{Linear and nonlinear hyperbolic systems}

\subsection*{Linearization involving one solution}

We begin with: 
 
\begin{theorem}[Existence theory for the Cauchy problem \cite{LeFloch-prepare}]  
\label{the} 
Given an entropy solution $ u$ (with small total variation) to a strictly hyperbolic 
system  
$$
 \del_t u + \del_x f(u) = 0
$$
and given an initial data $\psi_0 \in \Mcal_b(\RR)$, the initial value problem
associated with the linear hyperbolic system 
$$
\del_t \psi + \del_x (\widehat{Df}(u) \, \psi) = 0
$$
admits a (possibly non-unique) weak solution $\psi(t) \in \Mcal_b(\RR)$, 
satisfying on every compact time interval 
$$ 
\| \psi(t)\|_{\Mcal_b(\RR)} \lesssim \| \psi_0\|_{\Mcal_b(\RR)}. 
$$
\end{theorem}

For an earlier result in the genuinely nonlinear case, see Crasta and LeFloch \cite{CrastaLeFloch}. 

The system under consideration above contains compressive shocks, only.
Note also that solutions are understood in the sense of Volpert or, more generally, Dal~Maso-LeFloch-Murat \cite{DLM}. 
Our proof of Theorem~\ref{the} uses the property of pointwise convergence mentioned earlier
in Section~\ref{section22} and the following two lemmas.

\begin{lemma}[Continuous part]
There exists a time-dependent, bounded measure $ \mu= \mu(t)$, 
supported on the (countable) set of jump points of $ x \mapsto u(t,x)$, such that
for  
a.e. $ t \geq 0$ 
$$ 
\widehat{Df}(u^h) \, \psi^h 
\to  
\widehat{Df}(u) \, \psi + \mu 
$$ 
in the weak-star sense of measures.  

\end{lemma}

\begin{lemma}[Jump part]
If $ t \in \bigl(\underline T^h, \overline T^h \bigr) \mapsto y^h(t)$ 
is an approximate shock curve of a genuinely nonlinear family in $ u^h$, converging to 
some shock curve 
$ t \in \bigl(\underline T, \overline T \bigr) \mapsto y(t)$, then 
for a.e. $ t \in \bigl(\underline T, \overline T \bigr)$
$$  
\lim_{h \to 0} \psi^h \big(\big\{ y^h(t) \big\} \big) 
= \psi \big( \big\{ y(t)\big\} \big).   
$$

\end{lemma} 

A weaker result holds for non-genuinely nonlinear characteristic fields. 
We note that, for almost every time $t$ and at every jump point $x$ of the limiting function $v$, 
the left- and right-hand traces $v_\pm = v_\pm(t,x)$  
and $u_\pm= u_\pm(t,x)$ of the limiting solutions $u, v$ satisfy a relation 
$$
v_+ - v_- = \gamma \, (u_+ - u_-) \in \RN,  
$$
where $\gamma = \gamma(t,x)$ is a scalar.

The Dirac mass solutions found in \cite{LeFloch-IMA1} are now refered to in the literature 
as ``delta-shock waves''; see Zhang and Zheng \cite{ZhangZheng},  
Danilov and Shelkovich \cite{DanilovShelkovich}, 
and Mitrovic and Nedeljkov \cite{MN}. 
Furthermore, a large literature is available on solutions to linearized hyperbolic systems
involving a single entropy solution, especially by Bressan and Colombo \cite{BressanColombo}, 
Bressan, Crasta, and Piccoli \cite{Bressan1,Bressan2}, 
Ancona and Marson \cite{AnconaMarson2}, Bressan and Marson \cite{BressanMarson}, who investigated 
hyperbolic systems using the notion of ``first-order tangent vector'' associated with
piecewise smooth solutions; these papers encompass general entropy solutions to Temple class systems, 
as well as genuinely nonlinear systems. 


\subsection*{Linearization involving two solutions} 

The linearization of systems involving two solutions is comparatively more challenging. 
As far as the continuous dependence property of entropy solutions is concerned, 
one main result achieved by the linearization method is as follows: 

\begin{theorem}[Continuous dependence property \cite{HuLeFloch,LeFloch-book}]
\label{theo45} 
Consider the Cauchy problem for a strictly hyperbolic system of conservation laws with 
genuinely nonlinear fields. 
Then, any two entropy solutions $u,v$ with small total variation generated by front tracking 
depend $L^1$ continuously upon their initial data:  
$$
\|v(t) - u(t)\|_{L^1(\R)} \lesssim \|v(0) - u(0)\|_{L^1(\R)}.  
$$

\end{theorem}

Other proofs of the continuous dependence property were simultaneously and independently 
proposed by 
Liu and Yang \cite{LiuYang2} (Glimm scheme and explicit functional), 
and by Bressan, Liu, and Yang \cite{BressanLiuYang} (front tracking and explicit functional). 
The latter proofs are reviewed in the textbooks \cite{Bressan2,Dafermos2,HoldenRisebro}, 
while our proof via the linearization method included in the textbook \cite{LeFloch-book}
and was further developed in \cite{GoatinLeFloch2,LeFloch-cancel}. 
 
Our proof of Theorem~\ref{theo45}
is based a linearization approach and the construction of a weighted norm. 
A generalization to non-genuinely nonlinear systems of conservation laws is presented in \cite{LeFloch-prepare}
and uses the general existence theory established in \cite{IguchiLeFloch, GlassLeFloch}, 
as well as in 
\cite{Bianchini,BianchiniBressan}, \cite{LiuYang3}, \cite{AnconaMarson1,AnconaMarson3}. 

Let us sketch our proof of Theorem~\ref{theo45}. We begin by introducing the 
characteristic components $ \alpha = (\alpha_j)$ of a solution 
$$ 
\psi(t,x) =: \sum_j \alpha_j(t,x) \, \rt_j(t,x), \qquad (t,x) \in \Ccal(\At) \cap \Ccal(\psi),  
$$
and the weighted norm  
$$  
\| \psi(t)\|_{w(t)}  := \int_\RR \sum_j |\alpha_j(t,x)| \, w_j(t,x) \, dx,  
$$ 
where 
$ 0< w_{\min} \leq w_j(t,x) \leq w_{\max}$. For simplicity, we restrict the presentation here to the case that 
$\At = \Abar(u,v)$, the averaged matrix associated with two entropy solutions. 
  
It is not difficult to compute the time-derivative of the weighted norm. 
We suppose that the weight $w$ formally solves the adjoint system. Then, 
for any piecewise constant $\psi$ the weighted norm satisfies 
$$
\aligned 
 {d \over dt} \| \psi(t) \|_{w(t)}  
\leq & 
 \sum_{i,j}  \sum_{x \in \Jcal(\At(t))} \beta_j^-(t,x) \, w_j^-(t,x) + \beta_j^+(t,x) \, w_j^+(t,x) 
\endaligned 
$$
for all but finitely many $ t$, with 
$$
\aligned 
&  \beta_j^-(t,x) := \bigl(\lamb(t,x) - \lamt_{j-}(t,x) \bigr) 
  \, |\alpha_j^-(t,x)|,   
\\ 
&  \beta_j^+(t,x) := \bigl(\lamt_{j+}(t,x) - \lamb(t,x)\bigr) 
  \, |\alpha_j^+(t,x)|, 
\endaligned
$$
which we call the characteristic flux of the solution $ \psi$.

\begin{lemma}[Signs of the characteristic flux] 
Across each $ i$-shock, one has for $ j \neq i$ 
$$
\aligned 
&  \pm \beta_j^\pm \leq 0, \quad j<i, 
\\
&  \pm \beta_j^\pm \geq 0, \quad j>i, 
\endaligned 
$$
and for $ j= i$ 
$$
\aligned 
&      \hskip.35cm  \beta_i^\pm \leq 0,  \quad \Lcal_i, 
       \qquad \qquad 
       &  \beta_i^\pm \geq 0,  \quad \Rcal_i,   
\\
& \pm \beta_i^\pm \geq 0,  \quad \Scal_i,   
\qquad \qquad 
      &  \pm \beta_i^\pm \leq 0,  \quad  \Fcal_i. 
\endaligned
$$

\end{lemma}

This leads us to the following conditions on the weights: 

\begin{enumerate}

\item[$\bullet$] Case $j \neq i$: one always has one favorable sign and one unfavorable sign.  

\item[$\bullet$] Case $ j =i$: only rarefaction shocks correspond to two unfavorable signs !  

\end{enumerate}

\begin{definition} \it 
Fix a small $\kappa>0$. 
Consider a solution $\psi= \psi(t,x)$ together with its characteristic flux $\beta_j$. 
Then, at any $i$-shock the $j$-characteristic flux $(1 \leq j \leq N$) is said to be {\rm dominant} if 
$$
\kappa \, |\beta_j^-| \geq | \rt_{i+} - \rt_{i-}| \, |\beta_i^-| + | \At^+ - \At^-| \, \sum_k |\beta_k^-|. 
$$
\end{definition}

The dominant components enjoy the following properties.  

\begin{lemma} [Signs of the characteristic components]
For all $j \neq i$ 
$$
\sgn(\alpha_j^+) = \sgn(\alpha_j^-),     \quad \text{ $j$ dominant,}  
$$
while for $j=i$ 
$$
\sgn(\alpha_i^+) = 
\begin{cases} 
\hskip.25cm \sgn(\alpha_i^-), \quad & \text{ $\Lcal_i \cup \Rcal_i$ and $i$ dominant,}  
\\
-\sgn(\alpha_i^-), \quad & \text{ $\Scal_i \cup \Fcal_i$ and $i$ dominant.}  
\end{cases} 
$$
\end{lemma}

Hence, we see that the {\sl change of sign} of the characteristic component 
$\alpha$ is directly 
related to the nature (shock/rarefaction 
or undercompressive) of the discontinuity in $\At$. 

\begin{proposition}[Fundamental property nonlinear hyperbolic systems] 
Given two entropy solutions $ u,v$ with small total variation to 
a strictly hyperbolic, genuinely nonlinear system and let $ \Abar$ be the averaged matrix. 
Then, the matrix $ \At(t,x) := \Abar(u, v)(t,x)$
may contain compressive and undercompressive shocks 
but can not contain rarefaction-shocks. 

\end{proposition}

The proof based on a monotonicity property for the eigenvalue $ \lambda_i$:  
\begin{enumerate}

\item[--] For all $ u_-, v$  varying in a small neighborhood of $ 0$, 
the averaged speed $ \lamb_i(\cdot,v)$ is {\sl strictly monotone}
along the $ i$-shock curve from $ u_-$.  

\item[--] Moreover, if the $ i$-shock $ (u_-,u_+)$ satisfies Lax shock inequalities 
$$ 
\lam_i(u_-) >  \lamb_i(u_-, u_+) > \lam_i(u_+), 
$$
then the averaged speed satisfies 
$$ 
\lamb_i(u_-,v) > \lamb_i(u_+,v). 
$$
\end{enumerate}

Let us sketch the proof of these properties.  
The right-hand state $u_+=u_+(\eps)$ can be viewed as a function of $u_-$ and a parameter $\eps$
varying in the neighborhood of $0$, with 
$$ 
u_+(\eps) = u_- + \eps \, r_i(u_-) + O(\eps^2). 
$$ 
Then, we compute 
$$
\aligned 
\lamb_i(u_+,v) - \lamb_i(u_-,v)
& = \eps \, \nabla_1 \lamb_i(u_-, v) \cdot r_i(u_-) + O(\eps^2)
\\
& = \eps \, \nabla_1\lamb_i(u_-,u_-) \cdot r_i(u_-) + O(\eps^2) + O(\eps \, |v - u_-|) > 0. 
\endaligned  
$$
Using now that $\lamb_i(u,u) = \lam_i(u)$ and the symmetric property $\Abar(u,v) = \Abar(v,u)$, 
we obtain $2 \, \nabla_1 \lamb_i \cdot r_i = \nabla \lam_i \cdot r_i >0$.  
Provided $|\eps| + |v- u_-| \lesssim \delta_1$ is sufficiently small, 
we conclude that the function $\lamb_i(\cdot,v)$ is strictly monotone along the shock curve. 

If the shock $(u_-,u_+)$ satisfies the entropy inequalities 
and the normalization $\nabla\lam_i \cdot r_i >0$ is chosen, then $\eps < 0$ 
and we conclude from the above calculation that the averaged speed decreases from $u_-$ to $u_+$.

\section{Fluid dynamics equations} 

Finally, we turn the discussion to the compressible fluid equations with general equations of state, in 
either the form of the Lagrangian $p$-system  
\be
\label{36}
\del_t u_1 - \del_x u_2 = 0, \qquad \del_t u_2 + \del_x p(u_1) = 0. 
\ee
with $ p'(u_1) <0$, or equivalently in the form of the Euler equations 
\be
\label{37}
\del_t u_1 + \del_x (u_1 u_2) = 0, \qquad \del_t (u_1 u_2) + \del_x (u_1 (u_2)^2 + p(u_1)) = 0,  
\ee
with now $ p'(u_1) >0$.

\begin{theorem}[Continuous dependence property for fluid dynamics \cite{LeFloch-prepare}] 
\label{555} 
Consider solutions generated by front tracking and 
with small total variation of the fluid dynamics equations \eqref{36} and \eqref{37} 
Then, any two entropy solutions $u,v$ together with their approximations $ u_h, v_h$ 
satisfy for all $t \geq 0$ 
$$
\|v_h(t) - v_h(t)\|_{L^1(\R)} \lesssim \|v_h(0) - v_h(0)\|_{L^1(\R)} + o(h) 
$$
and, in consequence,  
$$
\|v(t) - u(t)\|_{L^1(\R)} \lesssim \|v(0) - u(0)\|_{L^1(\R)}.  
$$

\end{theorem}

More precisely, we prove the following sharp continuous dependence estimate 
$$  
\aligned
& \|v(t) - u(t) \|_{L^1(\RR)} + \Dbold_2(t) 
+ \Dbold_3(t) 
\lesssim \|v(0) - u(0)\|_{L^1(\RR)}, 
\endaligned
$$
where
$$
\Dbold_2(t) :=  \int_0^t \sum_{\Lcal(\at)}  | \lam^{\at} - \at_-| \, |v_- - u_-| \, d\tau,
$$
with a similar expression for $\Dbold_3(t)$. 
Similarly to the case of scalar equations, the terms $\Dbold_2(t)$ and 
$ \Dbold_3(t)$
are associated with compressive and undercompressive shocks in {\sl both} characteristic families, 
respectively.

The main steps of the proof of Theorem~\ref{555} are as follows:  
\begin{enumerate} 

\item The existence part follows from Iguchi-LeFloch's theory of non-genuinely nonlinear systems \cite{IguchiLeFloch}. 

\begin{enumerate} 

\item[--] One first approximates the pressure function by a function with finitely many inflection points. 

\item[--] One constructs interaction functionals adapted to the problem. 

\item[--] One derives uniform estimates that are independent of the number of inflection points. 

\end{enumerate}

\item The continuous dependence part relies on the following steps:  

\begin{enumerate} 

\item[--]  A classification of discontinuities based on the density/specific volume variable is 
established. 

\item[--] We show that the averaged matrix $ \Abar(u,v)$ contains no rarefaction shocks. 

\item[--] A weighted $ L^1$ norm is obtained by a constructive method by solving an evolution equation. 
\end{enumerate}

\end{enumerate}

The key observation of the linearization method to apply to fluid dynamics equations 
is provided by:

\begin{proposition}[Fundamental property for fluid dynamics] \it
If $ u, v$ are two entropy solutions to the fluid dynamics equations \eqref{36} or 
\eqref{37} with general equation of state for the pressure and for {\rm arbitrary large} total variation,  
then the averaged matrix $\At(t,x) := \Abar(u, v)(t,x)$ can not contain rarefaction-shocks.

\end{proposition}

\section{DLM theory of nonconservative hyperbolic systems} 
\label{DLMsection}
  
To complete this presentation, we want to stress the importance of small-scale phenomena for
formulating a well-posed hyperbolic theory. Recall the 
notion of family of paths in the sense of Dal~Maso-LeFloch-Murat \cite{DLM}. 
Let $ \Phi : [0,1] \times \RN \times \RN \times \RN$ be a map satisfying:  

\begin{enumerate}

\item $ \Phi(\cdot; u_-, u_+)$ is a path connecting $ u_-$ to $ u_+$ 
$$ 
\Phi(0; u_-, u_+) = u_-, \qquad \Phi(1; u_-, u_+) = u_+,
$$
$$ 
TV_{[0,1]}\big(\Phi (\cdot; u_-, u_+) \big) \lesssim |u_+ - u_-|;  
$$ 

\item $ \Phi$ is Lipschitz continuous in the graph distance 
$$  
\dist\Big( \Phi (\cdot; u_-, u_+), \Phi (\cdot; u_-', u_+') \Big) 
\lesssim |u_- - u_-'| + |u_+ - u_+'|.  
$$
\end{enumerate}

\begin{definition} \it
Given $ u \in BV(\RR, \RN)$ and 
$ g$ a Borel function, there exists a unique measure, called a {\rm nonconservative product}
and denoted by    
$$ 
\mu = \Big [g(u) \, \del_x u \Big ]_\Phi,
$$   
that is uniquely defined by the two properties: 

\begin{enumerate}
\item If $ B$ is a Borel subset of $ \Ccal(u)$, then
$ \mu(B) := \int_B g(u) \, \del_x u$. 
  
\item If $x \in \Jcal(u)$ is a point of jump then, with $ u_\pm := u_\pm(x)$, 
$$ 
\mu(\big\{ x \big\}) := \int_{[0,1]} g(\Phi(\cdot ; u_-, u_+)) \, \del_s \Phi (\cdot ; u_-, u_+).
$$
\end{enumerate}

\end{definition}

This definition provides us with a notion of weak solutions to, in particular,
nonconservative systems, once a family of paths is prescribed.
In the conservative case, our definition is independent of the paths and is 
consistent with the distributional definition 
$$
\Big [\nabla h(u) \, \del_x u \Big ]_\Phi = \del_x h(u).
$$

\begin{definition} \it 
A bounded variation function $u$ is a {\rm weak solution in the DLM sense}
 if for every test-function $\theta$
$$
- \int_{\RR \times \RR_+} \del_t \theta \, u  \, dxdt
+ \int_{\RR_+} \int_\RR \Big [A(u) \, \del_x u \Big ]_\Phi \, dt = 0.
$$ 

\end{definition}

The Riemann problem for systems of conservation laws was first solved by Lax \cite{Lax}. 
For nonconservative system we have the following generalization of Lax's theorem.  

\begin{theorem}[Riemann problem \cite{DLM}] 
Given a nonconservative, strictly hyperbolic, genuinely nonlinear system and a family of paths $\Phi$,  
the Riemann problem admits an entropy solution (in the DLM sense) 
satisfying Lax shock inequalities. 

\end{theorem}

Generalized Hugoniot jump relations for nonconservative systems read: 
$$  
- \lam^u \, (u_+ - u_-) + 
\int_0^1 A(\Phi(\cdot; u_-, u_+)) \, \del_s \Phi (\cdot; u_-, u_+) = 0.  
$$
Note that wave curves are only Lipschitz continuous at the origin, which is in contrast with Lax's standard $C^2$ 
regularity result. The regularity of the wave curves is also investigated
in Hayes and LeFloch \cite{HayesLeFloch}, 
Bianchini and Bressan \cite{Bianchini,BianchiniBressan},
Iguchi and LeFloch \cite{IguchiLeFloch}, and 
Liu and Yang \cite{LiuYang3}.

The initial value problem for nonconservative systems was solved by LeFloch and Liu. 

\begin{theorem}[Existence theory for the Cauchy problem \cite{LeFlochLiu}]  
Consider solutions constructed by the Glimm scheme for nonconservative, strictly hyperbolic, genuinely nonlinear system. 
Then, $u^h=u^h(t,x)$ have uniformly bounded total variation and
 converge to an entropy solution
$u$ in the DLM sense and for all but countably many times
\be
\label{start}
\Big [A(u^h) \, \del_x u^h \Big ]_\Phi(t) \rightharpoonup
 \Big [A(u) \, \del_x u \Big ]_\Phi(t).
\ee
\end{theorem} 

Observe that the `almost everywhere' convergence with respect to the Lebesgue measure is not strong enough
to derive the nonlinear stability property \eqref{start} and the pointwise convergence properties recalled at the end of 
Section ~\ref{section22} are necessary. It is interesting to observe similarities and differences
between in carrying out the theory of weak solutions for conservative systems or for nonconservative systems. 
In their work on the vanishing viscosity method, Bianchini and Bressan \cite{BianchiniBressan}
were able to cover {\sl both} conservative and nonconservative systems by the same arguments.  
As far as Glimm-type methods are concerned, 
a major distinction between conservative and nonconservative systems 
must be pointed out, when non-GNL are allowed. 
The following superposition property, essential in the theory of conservative systems,
does not extend to nonconservative systems: if three states $u_l, u_m, u_r$ satisfy the two relations 
$$
- \lamb \, (u_m - u_l) + f(u_m) - f(u_l) 
= - \lamb \, (u_r - u_m) + f(u_r) - f(u_m) = 0
$$
for some speed $\lamb$, 
then one also has  
$$
- \lamb \, (u_r - u_l) + f(u_r) - f(u_l) = 0. 
$$ 
For further material on nonconservative systems we refer to \cite{LeFloch-IMA2,LeFlochTzavaras0,LeFlochTzavaras,CLMP}.


\bibliographystyle{amsalpha}

\end{document}